\documentclass[11pt]{amsart}
\usepackage{amsmath}
\usepackage{amssymb}

\theoremstyle{plain}
\newtheorem{theorem}{\bf Theorem}[section]
\newtheorem{proposition}[theorem]{\bf Proposition}
\newtheorem{lemma}[theorem]{\bf Lemma}

\newtheorem{conjecture}[theorem]{\bf Conjecture}

\theoremstyle{definition}

\newtheorem{definition}[theorem]{\bf Definition}

 \DeclareMathOperator{\supp}{supp}

\numberwithin{equation}{section}

\begin{document}

\title[On short zero-sum subsequences
of zero-sum sequences]{On short zero-sum subsequences\\
of zero-sum sequences}

\address{Center for Combinatorics, Nankai University, LPMC-TJKLC, Tianjin
300071, P. R. China} \email{fys850820@163.com}

\address{Center for Combinatorics, Nankai University, LPMC-TJKLC, Tianjin
300071, P. R. China} \email{wdgao1963@yahoo.com.cn}

\address{Department of Mathematics, Tianjin Polytechnic University, Tianjin
300160, P. R. China} \email{gqwang1979@yahoo.com.cn}

\address{Center for Combinatorics, Nankai University, LPMC-TJKLC, Tianjin
300071, P. R. China} \email{zhongqinghai@yahoo.com.cn}

\address{Department of Mathematics, Dalian Maritime University, Dalian, 116026, P. R. China} \email{jjzhuang1979@yahoo.com.cn}

\author{Yushuang Fan, Weidong Gao, Guoqing Wang, Qinghai Zhong and Jujuan Zhuang}

\thanks{*Corresponding author: Guoqing Wang \\ E-mail:
gqwang1979@yahoo.com.cn}

%\author{
%Yushuang Fan, \ Weidong Gao,\ Guoqing Wang, \ Qinghai Zhong
%\bigskip
%\\
%Center for Combinatorics, LPMC-TJKLC\\
%Nankai University, Tianjin, 300071, P. R. China\\
%\\
%\\
%Jujuan Zhuang\\
%Department of Mathematics\\
%Dalian Maritime University, Dalian, 116026, P. R. China\\
%}

\date{}

\begin{abstract} Let $G$ be a finite abelian group, and let
$\eta(G)$ be the smallest integer $d$ such that every sequence over
$G$ of length at least $d$ contains a zero-sum subsequence $T$ with
length $|T|\in [1,\exp(G)]$. In this paper, we investigate the
question whether all non-cyclic finite abelian groups $G$  share
with the following property: There exists at least one integer $t\in
[\exp(G)+1,\eta(G)-1]$ such that every zero-sum sequence of length
exactly $t$ contains a zero-sum subsequence of length in
$[1,\exp(G)]$. Previous results showed that the groups $C_n^2$
($n\geq 3$) and $C_3^3$ have the property above. In this paper we
show that more groups including the groups $C_m\oplus C_n$ with
$3\leq m\mid n$, $C_{3^a5^b}^3$, $C_{3\times 2^a}^3$, $C_{3^a}^4$
and $C_{2^b}^r$ ($b\geq 2$) have this property. We also determine
all $t\in [\exp(G)+1, \eta(G)-1]$ with the property above for some
groups including the groups of rank two, and some special groups
with large exponent.
\end{abstract}

\subjclass[]{ }

\maketitle

\section{ Introduction }

Let $G$ be an additive  finite abelian group. We call a zero-sum
sequence $S$ over $G$  a {\sl short zero-sum sequence} if  $1\leq
|S|\leq \exp(G)$. Let $\eta(G)$ be the smallest integer $d$ such
that every sequence $S$ over $G$ of length $|S|\geq d$ contains a
short zero-sum subsequence. Let $D(G)$ be the Davenport constant of
$G$, i.e., the smallest integer $d$ such that every sequence over
$G$ of length at least $d$ contains a nonempty zero-sum subsequence.
Both $D(G)$ and $\eta(G)$ are classical invariants in combinatorial
number theory. For detail on terminology and notation we refer to
Section 2.

By the definition of $\eta(G)$ we know that for every integer $t\in
[1, \eta(G)-1]$, there is a sequence $S$ over $G$ of length exactly
$t$ such that $S$ contains no short zero-sum subsequence.  H.
Harborth \cite{Harborth} proved that  $\eta(C_3^3)=17$. In 1969,
Emde Boas and D. Kruyswijk \cite{EK69} observed the following
interesting fact when they investigated the problem of determining
of  $D(G)$.

\medskip
\noindent {\bf Theorem A.} {\sl Every zero-sum sequence over $C_3^3$
of length exactly 14 contains a short zero-sum subsequence.}

\medskip

Theorem A just asserts that every sequence over $C_3^3$ of length
$14$ containing no short zero-sum subsequence has sum nonzero.  In
1997, the second author of this paper obtained the following similar
result to Theorem A when he investigated the unique factorization
problem in zero-sum:

\medskip

\noindent{\bf Theorem B.} (\cite{Gao97})  {\sl Let $q=p^n$ with $p$
a prime, and let $G=C_q^2$. Then every zero-sum sequence $S$ over
$G$ of length $|S|\in [2q,3q-2]$ contains a short zero-sum
subsequence.}

\medskip
It is well known that $\eta(C_m^2)=3m-2$ for every $m\in
\mathbb{N}$. So, it is naturally to ask the following

\medskip
\noindent {\bf Open problem 1.}  {\sl Whether all non-cyclic finite
abelian groups $G$ share with the following property: There exists
at least one integer $t\in [\exp(G)+1,\eta(G)-1]$ such that every
zero-sum sequence over $G$ of length exactly $t$ contains a short
zero-sum subsequence.}

\medskip
%In this paper we shall investigate $\eta(S)$ first for general
%sequences $S$ and then for zero-sum sequences $S$.

\medskip Let us first make some easy observation for Open problem 1. Note that for every $t\in [1,D(G)]$
there exists a minimal zero-sum sequence over $G$ of length $t$. So,
to consider whether  $G$ has the property above we may assume that
$t\in [D(G)+1,\eta(G)-1]$. For convenience we introduce

\begin{definition} Let $G$ be a non-cyclic finite abelian group.
Define $C_0(G)=\{t: D(G)+1\leq t\leq \eta(G)-1,$ every zero-sum
sequence over  G  of length exactly t contains a short  zero-sum
subsequence $\}$.
\end{definition}

If $G=C_{2}\oplus C_{2m}$ then $D(G)+1=2m+2$ and  $\eta(G)-1=2m+1$.
Therefore, by the definition we have $C_0(C_2\oplus
C_{2m})=\emptyset$. We believe that the groups $C_2\oplus C_{2m}$
are the only non-cyclic groups which offer a negative answer to Open
problem 1.

\begin{conjecture} \label{conjecture1}  Let $G$ be a non-cyclic finite
abelian group. If $G\neq C_{2}\oplus C_{2m}$ then $C_0(G)\neq
\emptyset.$
\end{conjecture}

 From Theorem A and Theorem B we know that Conjecture
\ref{conjecture1} holds true for $G=C_3^3$ and $G=C_q^2$ with $q\geq
3$ is a prime power. In this paper we shall confirm Conjecture
\ref{conjecture1} for more groups by showing

\begin{theorem} \label{mainthm1}  If $G$ is one of the following
groups then $C_0(G)\neq \emptyset$.

\begin{enumerate}
\item $G=C_{n_1}\oplus C_{n_2}$ where $3\leq n_1\mid n_2$.

\item $G=C_{3^a5^b}^3$\ where $a\geq 1$ or $b\geq 2$.

\item $G=C_{3\times 2^a}^3$ where $a\geq 4$.

\item $G= C_{3^a}^4$ where $a\geq 1$.

\item $G=C_{2^a}^r$ where $3\leq r\leq a$, or $a=1$ and $r\geq 3$.

\item $G=C_{k}^{3}$ where $k=3^{n_1}5^{n_2}7^{n_3}11^{n_4}13^{n_5}$, $n_1\geq 1$, $n_3+n_4+n_5\geq 3$, and $n_1+n_2\geq
11+34(n_3+n_4+n_5)$.
\item $G=C_{p^nm}\oplus H$ where $p$ is a prime, $H$ is a finite abelian
$p$-group satisfying $|H|\geq 3$ and $p^n\geq D(H)$, and $p^nm\geq
3$.
\end{enumerate}

\end{theorem}

We also determine $C_0(G)$ completely for the following groups.

\begin{theorem} \label{mainthm2}  Let $G$ be a non-cyclic finite abelian group.
Then,
\begin{enumerate}
\item $C_{0}(G)=[D(G)+1, \eta(G)-1]$ if $r(G)=2$.

\item $C_{0}(G)=[D(G)+1, \eta(G)-1]$ if $G=C_{p^nm}\oplus H$ with $p$ a prime, $H$ a $p$-group and $p^n\geq D(H)$.

\item $C_{0}(C_{3}^{4})=\{\eta(C_{3}^{4})-2, \eta(C_{3}^{4})-1\}=\{37,38\}$.

\item $C_{0}(C_{2}^{r})=\{\eta(C_{2}^{r})-3, \eta(C_{2}^{r})-2\}$, where $r\geq 3$.
\end{enumerate}

\end{theorem}

The rest of this paper is organized as follows: In Section 2 we
introduce some notations and prove some preliminary results; In
Section 3  we shall derive some lower bounds on $\min \{C_{0}(G)\}$;
In Section 4 we study $C_0(G)$ with focus on the groups $C_3^r$; In
Section 5 and 6 we shall prove Theorem \ref{mainthm2} and Theorem
\ref{mainthm1}, respectively; and in the final Section 7 we give
some concluding remarks and some  open problems.

\section{ Notations and some preliminary results}
Our notations and terminologies are consistent with \cite{GG06} and
\cite{GH06}. We briefly gather some key notions and fix the
notations concerning sequences over finite abelian groups. Let
$\mathbb{Z}$ denote the set of integers.  Let $\mathbb{N}$ denote
the set of positive integers,  and
$\mathbb{N}_{0}=\mathbb{N}\cup\{0\}$.
 For real numbers $a, b\in \mathbb R$, we set $[a, b]=\{x \in \mathbb{Z} : a\leq x\leq b\}$.
Throughout this paper, all abelian groups will be written
additively, and for $n, r \in \mathbb N$, we denote by $C_n$ a
cyclic group with $n$ elements, and denote by $C_n^r$ the direct sum
of $r$ copies of $C_n$.

Let $G$ be a finite abelian group and $\exp(G)$ its exponent. A
sequence $S$ over $G$ will be written in the form
\[
S = g_1 \cdot \ldots \cdot g_{\ell} = \prod_{g \in G} g^{\mathsf v_g
(S)} \,, \quad \text{with} \ \mathsf v_g (S) \in \mathbb N_0 \
\text{for all} \ g \in G \,,
\]
and we call
\[
|S| = \ell \in \mathbb N_0 \quad  \text{the {\it length} \ and }
\quad \sigma (S) = \sum_{i=1}^{\ell} g_i = \sum_{g \in G} \mathsf
v_g (S)g \in G \quad  \text{the {\it sum} of} \ S \,.
\]
Let ${\rm supp}(S)=\{g\in G:  \mathsf v_g(S)>0\}.$ We call $S$ a
{\sl square free sequence} if $\mathsf v_g(S)\leq 1$ for every $g\in
G$. So, a square free sequence over $G$ is actually a subset of $G$.
A sequence $T$ over $G$ is called a {\sl subsequence} of $S$ if
$v_g(T)\leq v_g(S)$ for every $g\in G$, and denote by $T|S$.
   For every $r\in [1,\ell]$,  define
  $$
  \sum_{\leq r}(S)=\{\sigma(T): \ T\mid S,\ 1\leq |T|\leq r\}
  $$
  and define
  $$
  \sum(S)=\{\sigma(T): \ T\mid S,\ |T|\geq 1\}.
  $$

The sequence $S$ is called
\begin{itemize}
\item a {\sl zero-sum sequence} if $\sigma(S)=0$.

\item a {\sl short\ zero-sum\ sequence} over $G$ if it is a zero-sum sequence of length $|S|\in [1, \exp(G)]$.

\item a {\sl short  free} sequence over $G$ if $S$ contains no short
zero-sum subsequence.
\end{itemize}

So, a zero-sum sequence over $G$ which contains no short zero-sum
subsequence will be called a zero-sum short free sequence over $G$.

\medskip
For every element $g\in G$, we set
$g+S=(g+g_{1})\cdot\ldots\cdot(g+g_{\ell})$. Every map of abelian
groups $\varphi:G\rightarrow H$ extents to a map from the sequences
over $G$ to the sequences over $H$ by
$\varphi(S)=\varphi(g_{1})\cdot\ldots\cdot\varphi(g_{\ell})$. If
$\varphi$ is a homomorphism, then $\varphi(S)$ is a zero-sum
sequence if and only if $\sigma(S)\in \ker(\varphi)$.

\medskip
We shall study $C_0(G)$ by using the following property
which was first introduced and investigated by Emde Boas and
Kruyswijk \cite{EK69} in 1969 for the groups $C_p^2$ with $p$ a
prime,  and was investigated in 2007 for the groups $C_n^r$ by the
second author, Geroldinger and Schmid \cite{GGS07}.

\medskip
\noindent{\bf Property C:}  {\sl We say the group $C_{n}^{r}$ has
property C if $\eta(C_{n}^r)=c(n-1)+1$ for some positive integer
 $c$, and every short  free sequence $S$ over $C_n^r$ of length $|S|=c(n-1)$   has the
 form $S=\prod_{i=1}^{c}g_{i}^{n-1}$ where $g_1,\ldots,g_c$ are pairwise distinct elements of $C_n^r$.}

 \medskip
 It is conjectured that every group of the form $C_n^r$ has Property
 C(see \cite{GG06}, Section 7). The following three lemmas connect
 $C_0(G)$ with Property C.

\medskip

\begin{lemma} \label{Prog} Let $G=C_{n}^{r}$ with $\eta(G)=c(n-1)+1$
for some $c\in \mathbb{N}$. If $c\leq n$ and if $G$ has Property $C$
then $\eta(G)-1\in C_{0}(G)$.
\end{lemma}

\proof Let $S$ be a zero-sum sequence  over $G$ of length
$|S|=\eta(G)-1=c(n-1)$. Assume to the contrary that, $S$ contains no
short zero-sum subsequence. Since $G$ has Property $C$,
$S=\prod_{i=1}^{c}g_{i}^{n-1}$ with $g_{i}\in G$ for every $i\in [1,
c]$. Therefore, $(n-1)(g_{1}+g_{2}+\cdots
+g_{c})=\sigma(S)=0=n(g_{1}+g_{2}+\cdots +g_{c})$. It follows that
$g_{1}+g_{2}+\cdots +g_{c}=0$. But $g_{1}g_{2}\cdot\ldots\cdot
g_{c}\mid S$ and $|g_{1}g_{2}\cdot\ldots\cdot g_{c}|=c\leq n$, a
contradiction. Hence, $\eta(G)-1\in C_{0}(G)$.\qed

\begin{lemma}  \label{thmmain1}  Let $G$ be a finite abelian group, and let $H$ be a proper subgroup of $G$ with $\exp(G)=\exp(H)\exp(G/H)$.
Suppose that the following conditions hold.

{\rm (i)} \ $\eta(G)=(\eta(H)-1)\exp(G/H)+\eta(G/H)$;

{\rm (ii)} \ $G/H\cong C_n^r$ has Property C;

{\rm (iii)} \ There exist $t_1\in [1, \exp(G/H)-1]$ and  $t_2\in
\{1,2\}$ such that $t_2\leq t_1$ and such that $[\eta(G/H)-t_1,
\eta(G/H)-t_2]\subset C_0(G/H)$.

\noindent Then,
$$[\eta(G)-t_1,
\eta(G)-t_2]\subset C_0(G).$$

\end{lemma}

\proof To prove this lemma, we assume to the contrary that there
 is a  zero-sum short free sequence
$S$ over $G$ of length $\eta(G)-t$ for some  $t\in [t_2, t_1]$. Let
$\varphi$ be the natural homomorphism from  $G$ onto  $G/H$.

Note that
\begin{equation}\label{equation length of S}
|S|=\eta(G)-t=(\eta(H)-1)\exp(G/H)+(\eta(G/H)-t).
\end{equation}
This allows us to take an arbitrary decomposition of $S$
\begin{equation}\label{equation S=a decomposition}
S=\left(\prod\limits_{i=1}^{\eta(H)-1}S_i\right)\cdot S^{'}
\end{equation}
 with
\begin{equation}\label{equation S=a decomposition condI}
|S_i|\in[1,\exp(G/H)]
\end{equation} and
\begin{equation}\label{equation S=a decomposition condII}
\sigma(S_i)\in \ker \varphi=H
\end{equation}
for every $i\in [1,\eta(H)-1]$.

Combining \eqref{equation length of S}, \eqref{equation S=a
decomposition}, \eqref{equation S=a decomposition condI} and
\eqref{equation S=a decomposition condII} we infer that
\begin{equation}\label{equation |S'| geq -t1}
|S^{'}|\geq \eta(G/H)-t\geq \eta(G/H)-t_1
\end{equation}
 and
\begin{equation}\label{equation varphi(S')=0}
\sigma\bigl(\varphi(S^{\prime})\bigr)=0.
\end{equation}

\textbf{Claim A.}  $\varphi(S')$ contains no zero-sum subsequence of
length in $[1,\exp(G/H)]$.

\noindent{\sl Proof of Claim A.}  Assume to the contrary that, there
exists a subsequence $S_{\eta(H)}$ (say) of $S'$ of length
$|S_{\eta(H)}|\in [1,\exp(G/H)]$ such that $\sigma(S_{\eta(H)})\in
\ker \varphi=H$. It follows that the sequence
$U=\prod\limits_{i=1}^{\eta(H)}\sigma(S_i)$ contains a zero-sum
subsequence $W=\prod_{i\in I}\sigma(S_i)$ over $H$ with $I\subset
[1,\eta(H)]$ and $|W|=|I|\in [1,\exp(H)]$. Therefore, the sequence
$\prod_{i\in I}S_i$ is a zero-sum subsequence of $S$ over $G$ with
 $1\leq  |\prod_{i\in I}S_i|\leq |I|\exp(G/H)\leq
\exp(H)\exp(G/H)=\exp(G)$, a contradiction. This proves Claim A.

By \eqref{equation |S'| geq -t1}, \eqref{equation varphi(S')=0},
Claim A and Condition (iii), we conclude that $$t_2=2$$ and
\begin{equation}\label{equation exact length of S'}
|S^{'}|=\eta(G/H)-1.
\end{equation}
This together with Condition (ii) implies that
\begin{equation}\label{equation form of varphi(S')}
\varphi(S')=x_1^{n-1}\cdot\ldots\cdot x_c^{n-1}
\end{equation}
where $c=\tfrac{\eta(G/H)-1}{n-1}$
 and $x_1,\ldots,x_c$ are pairwise
distinct elements of the quotient group $G/H$. So, we just proved
that every decomposition of $S$ satisfying conditions (\ref{equation
S=a decomposition condI}) and (\ref{equation S=a decomposition
condII}) has the properties (2.5)-(2.8).

\medskip
Since $t\leq t_1\leq \exp(G/H)-1$, it follows from \eqref{equation
length of S}, \eqref{equation S=a decomposition condI} and
\eqref{equation exact length of S'} that $|S_i|\in [2,\exp(G/H)]$
for all $i\in [1,\eta(H)-1]$. Moreover, since $t\geq t_2=2$, it
follows that there exists $j\in [1,\eta(H)-1]$ such that $|S_j|\leq
\exp(G/H)-1$. Without loss of generality we assume that $$|S_1|\in
[2,\exp(G/H)-1].$$

Suppose that there exists $h\in {\rm supp}(\varphi(S_1))\cap {\rm
supp}(\varphi(S^{'}))$. By \eqref{equation form of varphi(S')}, we
have that the sequence $S_1\cdot S'$ contains a subsequence
$S_1^{\prime}$ with $\varphi(S_1^{\prime})=h^n$. Let
$S^{''}=S_1\cdot S'\cdot {S_1^{\prime}}^{-1}$. We get a
decomposition
$S=S_1^{\prime}\cdot\left(\prod\limits_{i=2}^{\eta(H)-1}S_i\right)\cdot
S^{''}$ satisfying (\ref{equation S=a decomposition condI}) and
(\ref{equation S=a decomposition condII}) . But
$|S^{''}|=|S_1|+|S'|-|{S_1^{\prime}}^{-1}| \leq (n-1)+
(\eta(G/H)-1)-n=\eta(G/H)-2$, a contradiction on (\ref{equation
exact length of S'}).  Therefore,
$${\rm supp}(\varphi(S_1))\cap {\rm supp}(\varphi(S^{'}))=\emptyset.$$

Take a term  $g\mid S_1$. Since $\varphi(g)\notin
{\rm supp}(\varphi(S'))$ and $|S'\cdot g|=\eta(G/H)$, it follows from
Claim A that $S'\cdot g$ contains a subsequence $S_1^{\prime}$ with
\begin{equation}\label{equation g divide S1'}
g\mid S_1^{\prime}
\end{equation}
and
\begin{equation}\label{equation |S1prime| leq}
|S_1^{\prime}|\leq \exp(G/H)
\end{equation}
and \begin{equation}\label{equation sigma(S1prime)}
\sigma(S_1^{\prime})\in \ker \varphi.
\end{equation}
Let $S^{''}=S_1\cdot S'\cdot {S_1^{\prime}}^{-1}$. By
\eqref{equation form of varphi(S')}, \eqref{equation g divide S1'},
\eqref{equation |S1prime| leq} and \eqref{equation sigma(S1prime)},
we conclude that $|{\rm supp}(\varphi(S''))|\geq c+1$, a
contradiction with \eqref{equation form of varphi(S')}. This proves
the lemma. \qed

\medskip
From Lemma \ref{thmmain1}, we immediately obtain the following

\begin{lemma}\label{Theorem general for C_0}  Let $r\in \mathbb{N}$, and let $G_1=C_{n_1}^r$, $G_2=C_{n_2}^r$ and $G=C_{n_1 n_2}^r$.
Suppose that the following conditions hold.

{\rm (i)} \
$\frac{\eta(G_1)-1}{n_1-1}=\frac{\eta(G_2)-1}{n_2-1}=\frac{\eta(G)-1}{n_1
n_2-1}=c$ for some $c\in \mathbb{N}$;

{\rm (ii)} \ $G_2$ has Property C;

{\rm (iii)} \ There exist $t_1\in [1, n_2-1]$, $t_2\in \{1,2\}$
such that $t_{2}\leq t_{1}$ and such that $[\eta(G_2)-t_1,
\eta(G_2)-t_2]\subset C_0(G_2)$.

\noindent Then,
$$[\eta(G)-t_1,
\eta(G)-t_2]\subset C_0(G).$$

\end{lemma}

\section{Some lower bounds on $\min C_0(G)$}

In this section we shall prove the following
\begin{proposition} \label{zz} Let $G=C_{n}^{r}$ with  $n\geq 3, r\geq 3$,  and let $\alpha_{r}\equiv -2^{r-1}({\rm mod}\ n)$ with $ \alpha_{r}\in [0,n-1]$.
Then,

\begin{enumerate}

\item $C_{0}(G)\subset [(2^r-1)(n-1)-\alpha_{r}+1, \eta(G)-1]$ if
$\alpha_{r}\neq 0$.

\item $C_{0}(G)\subset \{(2^r-1)(n-1)-n, (2^r-1)(n-1)-n+1\}$ if
$\alpha_{r}=0$.
\end{enumerate}

\end{proposition}

 Note that  $\alpha_r\neq 0$ if and only if  $n\neq
2^k$,  or $n=2^k$ and  $r-1<k$; and $\alpha_r=0$ if and only if
$n=2^k$ and $k\leq r-1$.

\medskip
For every $r\in \mathbb{N}$, let
$$G=C_{n}^{r}=<e_{1}>\oplus \cdots \oplus <e_{r}>$$
with $<e_i>=C_n$ for every $i\in [1,r]$, and let
$$S_{r}=\prod_{\{i_1, \ldots, i_k\} \subset [1,r]}(e_{i_{1}}+\cdots
+e_{i_{k}})^{n-1}$$ where $\{i_1,\ldots, i_k\}$ runs over all
nonempty subsets of $[1,r]$.  We can regard $C_n^r$ as a subgroup of
$C_n^{r+1}$ and therefore $S_{r+1}$ has the following decomposition
$$
S_{r+1}=S_r(S_r+e_{r+1})e_{r+1}^{n-1}.
$$

\medskip
Since the proof of Proposition \ref{zz} is somewhat long, we split
the proof into lemmas begin with the following easy one

\begin{lemma} \label{pre} $S_r$ is a short free sequence over $C_n^r$ of length $|S_{r}|=(2^r-1)(n-1)$ and of sum $\sigma(S_{r})=-2^{r-1}(e_{1}+\cdots
+e_{r})=\alpha_r(e_1+\cdots +e_r)$.
\end{lemma}

\proof Obviously. \qed

\medskip

\begin{lemma}\label{zerosum} Let $G=C_n^r$ with $r\geq 2$. Then for every   $m\in [1, n-1]$ and every $i\in [1, r]$, the sequence $S_{r}(e_{i}^m)^{-1}(me_{i})$
 contains no short zero-sum subsequence.
\end{lemma}

\proof Without loss of generality, we assume that $i=r$.

Assume to the contrary that $S_{r}(e_{r}^m)^{-1}(me_{r})$ contains a
short zero-sum subsequence $U$. Since $S_r$ contains no short
zero-sum subsequence we infer that  $me_r\mid U$. Therefore,
$U=(me_r)U_0 (U_1+e_r) e_r^k$ with $U_0\mid S_{r-1}$ and $U_1\mid
S_{r-1}$ and $k\in [0,n-1]$. It follows that $U_0U_1$ is zero-sum
and $1\leq |U_0U_1|\leq n-1$. Since every element in $\supp
(S_{r-1})$ occurs $n-1$ times in $S_{r-1}$, it follows from
$|U_0U_1|\leq n-1$ that $U_0U_1\mid S_{r-1}$. Therefore, $U_0U_1$ is
a short zero-sum subsequence of $S_{r-1}$, a contradiction with
Lemma \ref{pre}. \qed

\medskip
Let $A$ be a set of zero-sum sequences over $G$. Define $${\mathcal {L}
}(A)=\{|T|: T\in A\}.$$

In  this section below we shall frequently use the following easy
observation.

\begin{lemma} \label{construction} Let $G$ be a finite abelian
group, and let $a,b\in \mathbb{N}$ with $a\leq b$. If there exists
 a set $A$ of zero-sum short free sequences over $G$ such
that $[a, b]\subset {\mathcal L}(A)$, then $C_0(G)\cap [a,
b]=\emptyset$.
\end{lemma}

\proof It immediately follows from the definition of $C_0(G)$. \qed

\begin{lemma} \label{imp} Let $G=C_n^r$ with $n, r\geq 3$. Then,

\begin{enumerate}
\item $C_{0}(G)\cap [|S_r|-(3n-3)-\alpha_{r},
|S_r|-\alpha_{r}]=\emptyset$ if $\alpha_r\neq 0$.

\item $C_0(G)\cap [|S_r|-(3n-3), |S_r|-(n+1)]=\emptyset$ if
$\alpha_r=0$.
\end{enumerate}

\end{lemma}

\proof  Recall that $|S_r|=(2^r-1)(n-1)$. We split the proof into
three steps.

\medskip
\textbf{Step 1.} In this step we shall prove that $$C_0(G)\cap
[|S_r|-(3n-3)-\alpha_r, |S_r|-(n+1)-\alpha_r]=\emptyset$$ no matter
 $\alpha_r=0$ or not.

Let
$$A={\big\{}S_r{\big(}(e_1+\cdots+e_r)^{\alpha_r}W e_3^m{\big)}^{-1}(me_3): W\mid S_2, \sigma(W)=0, m\in [1, n-1]{\big\}}.$$
It follows from Lemma \ref{zerosum} that every sequence in $A$ is
zero-sum short free.

Since ${\mathcal L}(\{W: W\mid S_2, \sigma(W)=0\})=[n+1, 2n-1]$, we
conclude easily  that $${\mathcal L}(A)=[|S_r|-(3n-3)-\alpha_r,
|S_r|-(n+1)-\alpha_r].$$ Now the result follows from Lemma
\ref{construction} and Conclusion 2 follows.

\medskip
\textbf{Step 2.} We show that $C_0(G)\cap [|S_r|-(n+\alpha_{r}),
|S_r|-(r-1)\alpha_{r}]=\emptyset$ for $\alpha_r\neq 0$.

Let
\begin{align*}
A_1&={\big\{}S_r{\big(}(e_1+e_2)^{\alpha_r}e_3^{\alpha_r}\cdot
\ldots \cdot e_r^{\alpha_r}e_1^m{\big)}^{-1}(me_1): m\in [1,
n-1]{\big\}}
\end{align*}
and

 \begin{align*}
 A_2=\big \{S_r{\big(}(e_1+e_2)^{\alpha_r}(e_1+e_3)
e_3^{\alpha_r-1}e_4^{\alpha_r}\cdot \ldots \cdot e_r^{\alpha_r}
e_1^{n-1}{\big)}^{-1} \big\}.
\end{align*}
It is easy to see that every sequence in $A_1\cup A_2$ is zero-sum
short free  by Lemma \ref{zerosum} and Lemma \ref{pre}. Note that
\begin{align*}{\mathcal L}(A_1)\cup{\mathcal L}(A_2)&=[|S_r|-(r-1)\alpha_r-n+2, |S_r|-(r-1)\alpha_r]\cup\{|S_r|-(r-1)\alpha_r-n+1\}\\
&=[|S_r|-(r-1)\alpha_r-n+1, |S_r|-(r-1)\alpha_r].
\end{align*}
Since $r\geq 3$, we have $|S_r|-(r-1)\alpha_r-n+1\leq
|S_r|-(n+\alpha_r)$. Therefore, ${\mathcal L}(A_1\cup A_2)={\mathcal
L}(A_1)\cup{\mathcal L}(A_2)\supset [|S_r|-(n+\alpha_r),
|S_r|-(r-1)\alpha_r]$. Again the result follows from Lemma
\ref{construction}.

\medskip
\textbf{Step 3.} We prove $C_0(G)\cap [|S_r|-(r-1)\alpha_{r},
|S_r|-\alpha_{r}]=\emptyset$ for $\alpha_r\neq 0$.

Let
\begin{align*}
A={\big\{}&S_r{\big(}(e_1+\cdots +e_r)^{k_1}(e_1\cdot \ldots \cdot
e_r)^{k_2}(e_1+\cdots +e_{k_3})e_{k_3+1}\cdot \ldots \cdot
e_r{\big)}^{-1}: \\&k_1\in [0,\alpha_r-1],  k_2\in [0,\alpha_r-1]
,k_1+k_2=\alpha_r-1,k_3\in [1, r]{\big\}}.
\end{align*}
Then every sequence in $A$ is zero-sum short free by Lemma
\ref{zerosum} and by Lemma \ref{pre}, and
\begin{align*}
{\mathcal
L}(A)&=\{|S_r|-k_1-rk_2-1-(r-k_3): k_1+k_2=\alpha_r-1, k_2\in [0,
\alpha_r-1], k_3\in [1,r]\}\\&=\{|S_r|-\alpha_r-((r-1)k_2+(r-k_3)):
k_2\in [0, \alpha_r-1], k_3\in [1,r]\}\\&= [|S_r|-r\alpha_r,
|S_r|-\alpha_r].
\end{align*}
Now the result follows from Lemma
\ref{construction} and   the proof is completed. \qed

\medskip
We also need the following old easy results.
\begin{lemma} \label{extra} The following two conclusions hold.

\begin{enumerate}
\item If $n$ is a power of two then $\eta(C_n^r)=(2^r-1)(n-1)+1$.
\cite{Harborth}

\item $D(C_n^3)\geq 3n-2$. \cite{Olson}
\end{enumerate}
\end{lemma}

\medskip

\begin{lemma} \label{lem3.1} Let $n,r \in \mathbb{N}$ with $n\geq 3$
and $r\geq 3$.  If $\alpha_{r}\neq 0$ then $C_{0}(G)\subset
[(2^r-1)(n-1)-\alpha_{r}+1, \eta(G)-1]$.

\end{lemma}

\proof It suffices to show that $C_0(G)\cap [n+1,
|S_r|-\alpha_r]=\emptyset$.

We proceed  by induction on $r$. Suppose first that $r=3$.

By Lemma \ref{imp} and the definition of $C_0(C_n^3)$, we only need
to prove
 $$C_0(G)\cap [D(C_n^3)+1, |S_3|-(3n-3)-\alpha_3-1]=\emptyset.$$ By Lemma \ref{extra} we have $D(C_n^3)+1\geq
3n-1$. So,  it suffices to prove that $$C_0(G)\cap [3n-1,
|S_3|-(3n-3)-\alpha_3-1]=C_0(G)\cap [3n-1,
4n-4-\alpha_3-1]=\emptyset.$$

If $n=3$, then $[3n-1, 4n-4-\alpha_3-1]=\emptyset$ and the result
follows.

Now assume $n\geq 4$. It follows from $\alpha_3\neq 0$ that
$n\geq5$. Thus, $\alpha_3=n-4$ and $[3n-1,
4n-4-\alpha_3-1]=\{3n-1\}$.

Let $T=(e_1+e_2)^2(e_1+e_3)^{n-1}e_1^{n-1}e_2^{n-2}e_3$. Then $T$ is
zero-sum short free over $C_n^3$ of length $|T|=3n-1$. Now the
result follows from Lemma \ref{construction}.  This completes the
proof for $r=3$.

\medskip
Now assume that $r\geq 4$. By the induction hypothesis  there exists
a set $A_{r-1}$ of zero-sum short free sequences over $C_n^{r-1}$
such that $${\mathcal L}(A_{r-1})=[n+1,
|S_{r-1}|-\alpha_{r-1}].$$ Recall that $C_n^{r-1}\subset
C_n^r=C_n^{r-1} \oplus \langle e_r \rangle $. Let
\begin{align*}
A_r={\big\{}W_2(W_1+e_r)e_r^\ell: W_1\in A_{r-1}, W_2\in A_{r-1},
\ell\in [0, n-1], |W_1|+\ell \equiv 0 ({\rm mod}\ n) {\big\}}.
\end{align*}
Then every sequence in $A_r$ is zero-sum short free over $C_n^r$. It
is easy to see that
$${\mathcal L}(A_{r-1})\cup {\mathcal L}(A_r)\supset [n+1, 2|S_{r-1}|-2\alpha_{r-1}].$$
Note that
\begin{align*}2|S_{r-1}|-2\alpha_{r-1}& = |S_r|-(n-1)-2\alpha_{r-1}\\ &\geq
|S_r|-3(n-1).
\end{align*}
Therefore, $${\mathcal L}(A_{r-1})\cup{\mathcal L}(A_r)\supset[n+1,
|S_r|-3(n-1)].$$ Now the result follows from Lemma \ref{imp}. \qed

\begin{lemma} \label{lem3.2} Let $n,r,k\in \mathbb{N}$ with $k\geq
2, r\geq k+1$ and $n=2^k$, and let $G=C_n^r$. Then, $C_{0}(G)\subset
\{(2^r-1)(n-1)-n, (2^r-1)(n-1)-n+1\}$.
\end{lemma}

\proof  By Lemma \ref{extra} we have ,
$$|S_{r}|=(2^r-1)(n-1)=\eta(G)-1.$$

So, it suffices to show that $C_0(G)\cap([n+1, \eta(G)-(n+2)]\cup
[\eta(G)-n+1, \eta(G)-1])=\emptyset$.

Since $r\geq k+1$ we have
$$\sigma(S_{r})=0.$$

\textbf{Step 1.} We show $C_0(G)\cap[n+1, |S_r|-(n+1)]=\emptyset.$

We proceed by induction on $r$. Suppose first that $r=k+1$.

If $r=k+1=3$, we only need to prove $C_0(G)\cap [3n-1,
4n-5]=\emptyset$ by Lemma \ref{imp} and Lemma \ref{extra}. Let
\begin{align*}
A=&\{
(e_1+e_2+e_3)(e_1+e_2)^{n-1}(e_1+e_3)^{n-m}(e_2+e_3)e_1^{m}e_2^{n-1}e_3^{m-2}:
m\in [2, n-1]\} \cup
\\ &\{(e_1+e_2)^2(e_1+e_3)^{n-1}e_1^{n-1}e_2^{n-2}e_3\}.
\end{align*}
Then every sequence in $A$ is zero-sum short free and ${\mathcal
L}(A)=[3n-1,\ 4n-3]$ and  we are done.

If $r=k+1>3$, we have $\alpha_{r-1}\neq0$ and $r-1\geq3$, then by
Lemma \ref{lem3.1} there exists a set $A$ of zero-sum short free
sequences over $C_n^{r-1}$ such that ${\mathcal L} (A) \supset [n+1,
|S_{r-1}|-\alpha_{r-1}]$.

Let
\begin{align*}
B=A\cup{\big\{}W_2(W_1+e_r)e_r^\ell: W_1\in A, W_2\in A, \ell\in [0,
n-1], |W_1|+\ell \equiv 0 ({\rm mod}\ n) {\big\}}.
\end{align*}
Since
$$|S_{r-1}|-\alpha_{r-1}+|S_{r-1}|-\alpha_{r-1}+\alpha_{r-1}-1=|S_{r}|-3n/2,$$
we have ${\mathcal L}(B)\supset[n+1, |S_{r}|-3n/2]$. It follows from Lemma
\ref{imp} that $C_0(C_n^{r})\cap[n+1, |S_{r}|-(n+1)]=\emptyset$.

Now assume that $r>k+1$. By the induction hypothesis, we conclude
that there exists a set $A$ of zero-sum short free sequences
over $C_n^{r-1}$ such that ${\mathcal L} (A) \supset [n+1,
|S_{r-1}|-(n+1)]$.

Define  a set $B$ of zero-sum short free sequences over $C_n^r$ as
follows
\begin{align*}
B={\big\{}W_2(W_1+e_r)e_r^\ell: W_1\in A, W_2\in A, \ell\in [0,
n-1], |W_1|+\ell \equiv 0 ({\rm mod}\ n) {\big\}}.
\end{align*}
It is easy to see that
$${\mathcal L}(B)\supset [|S_{r-1}|-n,2|S_{r-1}|-2(n+1)]=[|S_{r-1}|-n, |S_{r}|-(3n+1)].$$

Let
\begin{align*}
C_1&=\{T:T\mid S_2, \sigma(T)=0\};\\
C_2&=\{(e_1+e_3)^{n-m}e_1^{m-1}e_2^{n-1}(e_1+e_2)e_3^{m}: m\in
[1, n-1]\};\\
C_3&=\{(e_1+e_2)^2(e_1+e_3)^{n-1}e_1^{n-1}e_2^{n-2}e_3\};\\
C_4&=\{(e_1+e_2+e_3)(e_1+e_2)^{n-1}(e_1+e_3)^{n-m}(e_2+e_3)e_1^{m}e_2^{n-1}e_3^{m-2}:
m\in [2, n-1]\}.
\end{align*}

Then every sequence in $ \cup_{i=1}^{4}C_i$ is zero-sum short free.
Clearly,
\begin{align*}
{\mathcal L}(C_1)&=[n+1, 2n-1];\\
{\mathcal L}(C_2)&=[2n, 3n-2];\\
{\mathcal L}(C_3)&=\{3n-1\};\\
{\mathcal L}(C_4)&=[3n, 4n-3].
\end{align*}
Let $$C = \cup_{i=1}^{4}C_i.$$ Then, $$ {\mathcal L}(C)\supset[n+1,
4n-3].$$

% Observe that
%$$[n+1,|S_{r-1}|-(n+1)]\subset
%[n+1,|S_{r-1}|-\alpha_{r-1}].$$
%By the induction hypothesis  there
%exists a set $A$ of zero-sum short free sequences over $C_n^{r-1}$
%such that ${\mathcal L} (A) \supset [n+1, |S_{r-1}|-(n+1)]$.

 Let
$$D=\{S_rT'^{-1}: T'\in C\}.$$ Then  every sequence in $D$ is zero-sum short
free, and
\begin{align*}
{\mathcal L}(D)&\supset[|S_r|-(4n-3), |S_r|-(n+1)]\\&\supset
[|S_r|-3n, |S_r|-(n+1)].
\end{align*}
This completes the proof of Step 1.

%If
%
%$C_0(C_n^{r_0-1})\cap[n+1, |S_{r_0}|-\alpha_{r_0}]=\emptyset$.
%So $C_0(G)\cap[n+1, |S_r|-(n+1)]=\emptyset.$
%
%Then,
%
% $$[n+1, |S_r|-(n+1)]=[n+1, 6n-8].$$
%
%Let $A_1$ and $A_2$ be the same as in the proof of 1. Let
%\begin{align*}
%A_3&=\{(e_1+e_2)^2(e_1+e_3)^{n-1}e_1^{n-1}e_2^{n-2}e_3\};\\
%A_4&=\{T:T=(e_1+e_2+e_3)(e_1+e_2)^{n-1}(e_1+e_3)^{n-m}(e_2+e_3)e_1^{m}e_2^{n-1}e_3^{m-2}, m\in [2, n-1]\}.
%\end{align*}
%Then every sequence in $ \cup_{i=1}^{4}A_i$ is zero-sum short free.
%Clearly,
%\begin{align*}
%{\mathcal L}(A_3)&=\{3n-1\};\\
%{\mathcal L}(A_4)&=[3n, 4n-3].
%\end{align*}
%It follows that $$C = \cup_{i=1}^{4}{\mathcal L}(A_i)\supset[n+1,
%4n-3].$$ This together with Lemma \ref{imp} gives that
%$C_0(G)\cap[n+1, 6n-8]=\emptyset.$

%To construct the required sequence $T$, by \eqref{equation sum of
%Sr=0} and \eqref{equation length of Sr=}, it is equivalent to show
%that for every $\ell\in [n+1, (2^r-1)(n-1)-(n+1)]$, there exists a
%zero-sum subsequence $T^{'}$ of $S_{r}$ with $|T^{'}|=\ell$.

\medskip
\textbf{Step 2.} We prove $C_0(G)\cap [\eta(G)-n+1,
\eta(G)-1]=\emptyset$.

Let $$A=\{S_r(e_r^m)^{-1}(me_r): m\in [1, n-1]\}.$$ Then every
sequence in $ A$ is zero-sum short free and
$${\mathcal L}(A)=[|S_r|-n+2, |S_r|]=[\eta(G)-n+1, \eta(G)-1].$$
This completes the proof. \qed

\bigskip
\noindent{\sl Proof of Proposition \ref{zz}}. 1. It is just Lemma
\ref{lem3.1}.

\medskip
2. Since $\alpha_{r}=0$, we have $n=2^k$ for some $k\in [2,r-1]$,
now the result follows from Lemma \ref{lem3.2}. \qed

\section{On the groups $C_3^r$}
In this section we shall study $C_0(G)$ with focus on $G=C_3^r$. Let
us first introduce some notations.

\begin{definition}  Let $G$ be a finite abelian group. Define
$g(G)$  to be the smallest integer $t$ such that every square free
sequence over $G$ of length  $t$ contains a zero-sum subsequence of
length  $\exp(G)$. Let $f(G)$  be the smallest  integer $t$ such
that every square free sequence over $G$ of length $t$ contains a
short zero-sum subsequence.
\end{definition}

We now gather some known results on Property C, $\eta(G), g(G)$ and
$f(G)$ which will be used in the sequel.

\begin{lemma}\label{Lemma list of contants}  Let $r,t\in
\mathbb{N}$, and let $n\geq 3$ be an odd integer. Then,

\begin{enumerate}
\item $\eta(C_n^3)\geq 8n-7$. \cite{Lower bounds}

\item $\eta(C_n^4)\geq 19n-18$. \cite{affine caps}

\item $\eta(C_3^3)=17=8\times 3-7$. \cite{Kemnitz}

\item $\eta(C_3^4)=39=19\times 3-18$. \cite{Kemnitz}

\item $\eta(C_5^3)=33=8\times 5-7$. \cite{GHST}

\item $\eta(C_{2^t}^r)=(2^r-1)(2^t-1)+1$. \cite{Harborth}

\item $\eta(C_{3\times 2^\alpha}^{3})=7(3\times 2^\alpha-1)+1$ where $\alpha\geq 1$. \cite{GHST}

\item $C_5^3$ and $C_{2^t}^{r}$ has Property $C$.  \cite{GHST}

\item $\eta(C_3^r)=2f(C_3^r)-1$. \cite {Harborth}

\item   $C_3^r$ has Property C. \cite {Harborth}
\end{enumerate}
\end{lemma}

\begin{proposition}\label{coraa} Let $r, t \in \mathbb{N}$. Then,
\begin{enumerate}
\item $C_{0}(C_{3}^{3})\subset [\eta(C_{3}^{3})-4,
\eta(C_{3}^{3})-1]$.

\item $C_{0}(C_{5}^{3})\subset [\eta(C_{5}^{3})-5, \eta(C_{5}^{3})-1]$.

\item $\begin{array}{llll}C_{0}(C_{2^t}^{r})\subset
\left\{\begin{array}{llll}
{[\eta(C_{2^t}^{r})-(2^t-2^{r-1}), \eta(C_{2^t}^{r})-1]}, &\mbox{if \ \ }r\leq t,\\
{[\eta(C_{2^t}^{r})-(2^t+1), \eta(C_{2^t}^{r})-2^t]}, &\mbox{if \ \
}r>t.\\
\end{array}
\right.
\end{array}$

\item $C_{0}(C_{6}^{3})\subset \{34, 35\}$.
\end{enumerate}
\end{proposition}

\proof  Conclusions 1, 2 and 4 follow from Proposition \ref{zz}. So,
it remains to prove Conclusion 3. If $r\leq t$ then applying
Proposition \ref{zz} with $\alpha_{r}=2^t-2^{r-1}$ we get,
$C_{0}(C_{2^t}^{r})\subset [(2^r-1)(2^t-1)-(2^t-2^{r-1})+1,
\eta(C_{2^t}^{r})-1] =[\eta(C_{2^t}^{r})-(2^t-2^{r-1}),
\eta(C_{2^t}^{r})-1]$. If $r>t$ then applying Proposition \ref{zz}
with  $\alpha_{r}=0$ we get,
$C_{0}(C_{2^t}^{r})\subset[\eta(C_{2^t}^{r})-(2^t+1),
\eta(C_{2^t}^{r})-2^t]$.  \qed

\begin{lemma} \label{Lemma the uniqueness of affine cap} (\cite{affine caps}, Lemma 5.4)  Let $r\in [3,5]$,
let $S$ and $S'$ be two  square free sequences over $C_3^r$ of
length $|S|=|S'|=g(C_3^r)-1$. Suppose that both $S$ and $S'$ contain
no zero-sum subsequence of length $3$. Then  $S'=\varphi(S)+a$,
where $\varphi$ is an automorphism of $C_{3}^{r}$ and $a\in
C_{3}^{r}$.
\end{lemma}

\begin{lemma}\label{Lemma the sequence in C33 }   (\cite{Davenport}, Lemma
1)  Let $T$ be a square free sequence over $C_3^3$ of length 8. If
$T$ contains no short zero-sum subsequence then there exists an
automorphism $\varphi$ of $C_3^3$ such that $\varphi(T)
=\left(\!\!\!\begin{array}{c} 0 \\ 1 \\0\end{array}\!\!\!\right)
\left(\!\!\!\begin{array}{c} 0 \\ 0 \\ 1\end{array}\!\!\!\right)
\left(\!\!\!\begin{array}{c} 0 \\ 1 \\ 1\end{array}\!\!\!\right)
\left(\!\!\!\begin{array}{c} 1 \\ 0 \\ 0\end{array}\!\!\!\right)
\left(\!\!\!\begin{array}{c} 1 \\ 2 \\ 0\end{array}\!\!\!\right)
\left(\!\!\!\begin{array}{c} 1 \\ 1 \\ 1\end{array}\!\!\!\right)
\left(\!\!\!\begin{array}{c} 1 \\ 1 \\ 2\end{array}\!\!\!\right)
\left(\!\!\!\begin{array}{c} 2 \\ 0 \\ 1\end{array}\!\!\!\right)$.

\end{lemma}

\begin{lemma}\label{Lemma the sequence in C34}  (\cite{Game}; \cite{affine caps}, page 182)
The following square free sequence over $C_3^4$ of length $20$
contains no zero-sum subsequence of length $3$.

$\left(\!\!\!\begin{array}{c} 0 \\ 0 \\ 0 \\ 0\end{array}\!\!\!\right)
\left(\!\!\!\begin{array}{c} 2 \\ 0 \\ 0 \\ 0\end{array}\!\!\!\right)
\left(\!\!\!\begin{array}{c} 0 \\ 2 \\ 0 \\ 0\end{array}\!\!\!\right)
\left(\!\!\!\begin{array}{c} 2 \\ 2 \\ 0 \\ 0\end{array}\!\!\!\right)
\left(\!\!\!\begin{array}{c} 1 \\ 0 \\ 2 \\ 0\end{array}\!\!\!\right)
\left(\!\!\!\begin{array}{c} 0 \\ 1 \\ 2 \\ 0\end{array}\!\!\!\right)
\left(\!\!\!\begin{array}{c} 1 \\ 2 \\ 2 \\ 0\end{array}\!\!\!\right)
\left(\!\!\!\begin{array}{c} 2 \\ 1 \\ 2 \\ 0\end{array}\!\!\!\right)
\left(\!\!\!\begin{array}{c} 1 \\ 1 \\ 1 \\ 0\end{array}\!\!\!\right)
\left(\!\!\!\begin{array}{c} 1 \\ 1 \\ 0 \\ 1\end{array}\!\!\!\right)$

$\left(\!\!\!\begin{array}{c} 0 \\ 0 \\ 2 \\ 2\end{array}\!\!\!\right)
\left(\!\!\!\begin{array}{c} 2 \\ 0 \\ 2 \\ 2\end{array}\!\!\!\right)
\left(\!\!\!\begin{array}{c} 0 \\ 2 \\ 2 \\ 2\end{array}\!\!\!\right)
\left(\!\!\!\begin{array}{c} 2 \\ 2 \\ 2 \\ 2\end{array}\!\!\!\right)
\left(\!\!\!\begin{array}{c} 1 \\ 0 \\ 0 \\ 2\end{array}\!\!\!\right)
\left(\!\!\!\begin{array}{c} 0 \\ 1 \\ 0 \\ 2\end{array}\!\!\!\right)
\left(\!\!\!\begin{array}{c} 1 \\ 2 \\ 0 \\ 2\end{array}\!\!\!\right)
\left(\!\!\!\begin{array}{c} 2 \\ 1 \\ 0 \\ 2\end{array}\!\!\!\right)
\left(\!\!\!\begin{array}{c} 1 \\ 1 \\ 1 \\ 2\end{array}\!\!\!\right)
\left(\!\!\!\begin{array}{c} 1 \\ 1 \\ 2 \\ 1\end{array}\!\!\!\right)$.

\end{lemma}

\begin{lemma}\label{Lemma S containing no short zero}   Let $G=C_3^r$ with $r\geq 3$, and let $S$ be
a sequence over $G$. Then,
\begin{enumerate}
\item If $S$ is a short free sequence over $G$ of length
$|S|=\eta(G)-1$, then $\sum\limits_{\leq 2}(S)=C_3^r\setminus\{0\}.$

\item Let  $T$ be a square free and short  free
sequence over $G$, and let $S=T^2$.  Then, for every  $g\in {\rm
supp}(S)$ we have, $\sum\limits_{\leq 2}(S\cdot
g^{-1})=\sum\limits_{\leq 2}(S)\setminus \{2g\}$.

\item If every short  free sequence of length $\eta(G)-1$
 has sum zero, then
$\eta(G)-2\in C_0(G)$.
\end{enumerate}

\end{lemma}

\proof Conclusions 1 and 2 are obvious.

To prove Conclusion 3, we  assume to the contrary  that
$\eta(G)-2\notin C_0(G)$, i.e., there exists a zero-sum short free
sequence $S$ over $G$ of length $|S|=\eta(G)-2$. By Lemma \ref{Lemma
list of contants}, we have $\eta(G)-2=2(f(G)-2)+1$. This forces that
$S=g_1^2\cdot \ldots \cdot g_{f(G)-2}^2 g_{f(G)-1}$ for some
distinct elements $g_1, \ldots, g_{f(G)-1}$ with $g_1\cdot \ldots
\cdot g_{f(G)-1}$ contains no short zero-sum subsequence. Put
$T=S\cdot g_{f(G)-1}$. Then $|T|=\eta(G)-1$. But $T$ contains no
short zero-sum subsequence and $\sigma(T)=g_{f(G)-1}\neq 0$, a
contradiction. \qed

\begin{lemma} \label{Lemma every sequence in C33 has sum zero}  Every short  free sequence over $C_3^3$ of length $16$  has sum zero. \end{lemma}

\proof Let $S$ be an arbitrary short  free sequence  over $C_3^3$ of
length $|S|=16$. From Lemma \ref{Lemma list of contants} we obtain
that $S=T^2$, where $T$ is a square free and short  free sequence
over $C_3^3$ of length $8$. It follows from Lemma \ref{Lemma the
sequence in C33 }  that $\sigma(T)=0$. Therefore,
$\sigma(S)=2\sigma(T)=0$. \qed

\begin{lemma} \label{Lemma C33-2 and C34-2}   The following two conclusions
hold.

\begin{enumerate}
\item $\{14,15\}=\{\eta(C_3^3)-3,\eta(C_3^3)-2\}\subset C_0(C_3^3)$.

\item $\{37,38\}=\{\eta(C_3^4)-2,\eta(C_3^4)-1\}\subset C_0(C_3^4)$.
\end{enumerate}

\end{lemma}

\proof \noindent{1.} From Theorem A we know that $14\in C_0(C_3^3)$, and $15\in
C_0(C_3^3)$  follows from Lemma \ref{Lemma list of contants}, Lemma
\ref{Lemma S containing no short zero} and Lemma \ref{Lemma every
sequence in C33 has sum zero}.

\medskip
\noindent{2.} Denote by $U$  the square free sequence over $C_3^4$
given in Lemma \ref{Lemma the sequence in C34}. It follows from
Conclusion 4 of Lemma \ref{Lemma list of contants} that $U$ is a
square free sequence of maximum length which contains no zero-sum
subsequence of length $3$.

Choose an arbitrary square free sequence  $T$ over  $C_3^4$ of
length $f(C_3^4)-1=19$ such that $T$ contains no short zero-sum
subsequence.

\textbf{Claim A.}  $\sigma(T)\notin -{\rm supp}(T)\cup\{0\}$.

\noindent{\sl Proof of Claim A.}  Put $S=T\cdot 0$. It follows from
Conclusion 4 of Lemma \ref{Lemma list of contants}  that $S$ is a
square free sequence over $C_3^4$ of maximum length which contains
no zero-sum subsequence of length $3$. By Lemma \ref{Lemma the
uniqueness of affine cap}, there exists an automorphism $\varphi$ of
$C_3^4$ and some $g\in C_3^4$ such that $S=\varphi(U-g)$. Since
$0\mid S$, it follows that $g\mid U$. Thus,
$\sigma(T)=\sigma(S)=\sigma(\varphi(U-g))=\varphi(\sigma(U-g))=\varphi(\sigma(U)-20g)=\varphi(\sigma(U)+g)$.
It is easy to check that $\sigma(U)={\tiny\left(\!\!\!\begin{array}{c} 2 \\ 2 \\ 2 \\
2\end{array}\!\!\!\right)}$. Since $-\sigma(U)={\tiny\left(\!\!\!\begin{array}{c} 1 \\
1 \\ 1 \\ 1\end{array}\!\!\!\right)}\notin {\rm supp}(U)$, it follows
that $-\sigma(T)=-\varphi(\sigma(U)+g)=\varphi(-\sigma(U)-g)\notin
\varphi({\rm supp}(U)-g)={\rm supp}(S)={\rm supp}(T)\cup\{0\}$. This
proves Claim A.

\medskip
From Lemma \ref{Lemma list of contants} and Claim A, we
derive that every short  free sequence over $C_3^4$ of length
$\eta(C_3^4)-1=38$ has sum nonzero. This is equivalent to that every
zero-sum sequence over $C_3^4$ of length $\eta(C_3^4)-1$ contains a
short zero-sum subsequence. Hence, $38=\eta(C_3^4)-1\in C_0(C_3^4).$

\medskip
Suppose that $37=\eta(C_3^4)-2\notin C_0(G)$, that is, there exists
a zero-sum short free sequence $V$ over $C_3^4$ of length
$|V|=\eta(C_3^4)-2=37$. From Lemma \ref{Lemma list of contants} we
know that $V=W^2 h^{-1}$, where $h\mid W$ and $W$ is a square free
and short free sequence over $G$ of length $f(C_3^4)-1=19$. It
follows from $\sigma (V)=0$ that $\sigma(W)=-h\in -{\rm supp}(W)$, a
contradiction with Claim A.  \qed

\begin{proposition}\label{Prop cap C_0(G)}   Let $G=C_3^r$ with $r\geq 3$.
If there is a short free sequence $S$ over $G$ of length
$|S|=\eta(G)-1$ such that $\sigma(S)\neq 0$, then
\begin{enumerate}
\item $|\{\eta(G)-2,\eta(G)-3\}\cap C_0(G)|\leq 1$

\item $|\{\eta(G)-3,\eta(G)-4\}\cap C_0(G)|\leq 1$.
\end{enumerate}
\end{proposition}

\proof \noindent{1.} Since $\sigma(S)\neq 0$, it follows from Lemma \ref{Lemma S
containing no short zero} that there exists a subsequence $W$ of $S$
of length $|W|\in \{1,2\}$ such that $\sigma(S)=\sigma(W)$.
Therefore, $\sigma(S\cdot W^{-1})=0$, $|S\cdot W^{-1}|\in
\{\eta(G)-3,\eta(G)-2\}$ and $S\cdot W^{-1}$ contains no short
zero-sum subsequence. Hence, $\eta(G)-2\notin C_0(G)$ or
$\eta(G)-3\notin C_0(G)$.

\medskip
\noindent{2.} By Lemma \ref{Lemma list of contants}, we have that
$S=T^2$, where $T$ is a square free sequence over $G$. Choose $g\in
{\rm supp}(S)$ such that $\sigma(S\cdot g^{-1})\neq 0$. Since
$\sigma(S\cdot g^{-1})=\sigma(S)-g\neq 2g$, it follows from Lemma
\ref{Lemma S containing no short zero} that $\sigma(S\cdot
g^{-1})\in \sum\limits_{\leq 2}(S\cdot
g^{-1})=C_3^r\setminus\{0,2g\}$. Similarly to Conclusion 1, we infer
that $\eta(G)-3\notin C_0(G)$ or $\eta(G)-4\notin C_0(G)$. \qed

\begin{proposition} \label{propaa}  $C_0(C_3^4)=\{37,38\}$.
\end{proposition}

\proof   By Proposition \ref{zz}, we have
\begin{equation}\label{equation C0(G) cap no [, 29]}
C_{0}(C_{3}^{4})\subset [30, \eta(C_3^4)-1]=[30,38].
\end{equation}
We show next that
\begin{equation}\label{equation C0(G) cap no [30,35]}
[30,36]\cap C_0(C_3^4)=\emptyset.
\end{equation}
Put
\begin{align*}
&T_2=\left(\!\!\!\begin{array}{c}
2 \\ 2 \\ 2 \\ 2\end{array}\!\!\!\right)^2;\\
&T_3=\left(\!\!\!\begin{array}{c}
2 \\ 2 \\ 0 \\ 0\end{array}\!\!\!\right) \left(\!\!\!\begin{array}{c} 0 \\ 0 \\
2 \\ 2\end{array}\!\!\!\right) \left(\!\!\!\begin{array}{c} 2 \\ 2 \\ 2 \\
2\end{array}\!\!\!\right);\\
&T_4=\left(\!\!\!\begin{array}{c} 2 \\ 0 \\ 0 \\
0\end{array}\!\!\!\right) \left(\!\!\!\begin{array}{c} 0 \\ 2 \\ 0 \\
0\end{array}\!\!\!\right) \left(\!\!\!\begin{array}{c} 0 \\ 0 \\ 2 \\
2\end{array}\!\!\!\right) \left(\!\!\!\begin{array}{c} 2 \\ 2 \\ 2 \\
2\end{array}\!\!\!\right);\\
&T_5=\left(\!\!\!\begin{array}{c} 2 \\ 0 \\ 0 \\
0\end{array}\!\!\!\right)^2 \left(\!\!\!\begin{array}{c} 0 \\ 2 \\ 0 \\
0\end{array}\!\!\!\right) \left(\!\!\!\begin{array}{c} 0 \\ 0 \\ 2 \\
2\end{array}\!\!\!\right) \left(\!\!\!\begin{array}{c} 0 \\ 2 \\ 2 \\
2\end{array}\!\!\!\right);\\
&T_6=\left(\!\!\!\begin{array}{c} 0 \\ 2 \\ 0 \\
0\end{array}\!\!\!\right) \left(\!\!\!\begin{array}{c} 2 \\ 2 \\ 2 \\
2\end{array}\!\!\!\right)^2 \left(\!\!\!\begin{array}{c} 0 \\ 1 \\ 0 \\
2\end{array}\!\!\!\right) \left(\!\!\!\begin{array}{c} 1 \\ 2 \\ 0 \\
2\end{array}\!\!\!\right) \left(\!\!\!\begin{array}{c} 2 \\ 1 \\ 0 \\
2\end{array}\!\!\!\right);\\
&T_7=\left(\!\!\!\begin{array}{c} 2 \\ 0 \\ 0 \\
0\end{array}\!\!\!\right) \left(\!\!\!\begin{array}{c} 0 \\ 2 \\ 0 \\
0\end{array}\!\!\!\right) \left(\!\!\!\begin{array}{c} 1 \\ 0 \\ 2 \\
0\end{array}\!\!\!\right) \left(\!\!\!\begin{array}{c} 0 \\ 1 \\ 2 \\
0\end{array}\!\!\!\right) \left(\!\!\!\begin{array}{c} 1 \\ 2 \\ 2 \\
0\end{array}\!\!\!\right) \left(\!\!\!\begin{array}{c} 0 \\ 0 \\ 2 \\
2\end{array}\!\!\!\right) \left(\!\!\!\begin{array}{c} 0 \\ 2 \\ 2 \\
2\end{array}\!\!\!\right);\\
&T_8=\left(\!\!\!\begin{array}{c} 2 \\ 0 \\ 0 \\
0\end{array}\!\!\!\right) \left(\!\!\!\begin{array}{c} 0 \\ 2 \\ 0 \\
0\end{array}\!\!\!\right)^2 \left(\!\!\!\begin{array}{c} 1 \\ 0 \\ 2 \\
0\end{array}\!\!\!\right) \left(\!\!\!\begin{array}{c} 0 \\ 1 \\ 2 \\
0\end{array}\!\!\!\right) \left(\!\!\!\begin{array}{c} 1 \\ 2 \\ 2 \\
0\end{array}\!\!\!\right) \left(\!\!\!\begin{array}{c} 0 \\ 0 \\ 2 \\
2\end{array}\!\!\!\right)^2.
\end{align*}
Let $U$ be the square free sequence given in Lemma \ref{Lemma the
sequence in C34}. Then $\sigma(U)=\left(\!\!\!\begin{array}{c} 2 \\ 2
\\ 2 \\ 2\end{array}\!\!\!\right)$. Let $S=U^2\cdot 0^{-2}$. We see
that $S$ is a short  free sequence of length $38=\eta(C_3^4)-1$. By
removing $T_i$ from $S$, we obtain that the resulting sequence $S_i$
is a zero-sum short free sequence of length $\eta(G)-i-1=38-i$. This
proves \eqref{equation C0(G) cap no [30,35]}. Combining
\eqref{equation C0(G) cap no [, 29]}, \eqref{equation C0(G) cap no
[30,35]} and Lemma \ref{Lemma C33-2 and C34-2}. we conclude that
$C_0(C_3^4)=\{\eta(G)-2,\eta(G)-1\}=\{37,38\}$. \qed

\section{Proof of Theorem \ref{mainthm2}}

In this section we shall prove Theorem \ref{mainthm2} and we need
the following lemmas.
\begin{lemma} \label{nd}(\cite{GGS}, \mbox{Theorem 5.2})  Every sequence $S$ over $C_{n}^{2}$ of length  $|S|=3n-2$ contains a zero-sum subsequence of length $n$ or $2n$.
\end{lemma}

\begin{lemma}\label{Lemma pho(G) and pho(H)} (\cite{GGS07}, \mbox{Lemma 4.5})  Let $G$ be a finite abelian group, and let $H$ be a proper subgroup of $G$ with
 $\exp(G)=\exp(H)\exp(G/H)$.
Then $\eta(G)\leq (\eta(H)-1)\exp(G/H)+\eta(G/H)$.
\end{lemma}

\begin{lemma} \label{Davenportconstant}  Let $p$ be a prime and let $H$ be a finite abelian $p$-group such that $p^n\geq D(H)$.
Let $n_{1}, n_{2}, m, n\in \mathbb{N}$ with $n_{1}\mid n_{2}$. Then,
\begin{enumerate}
\item $D(C_{n_{1}}\oplus C_{n_{2}})=n_{1}+n_{2}-1$. (\cite{Olson})

\item $D(C_{mp^n}\oplus H)=mp^n+D(H)-1$. (\cite{EK69})

\item Let $G=C_{p^{e_1}}\oplus \cdots \oplus C_{p^{e_r}}$ with $e_i\in \mathbb{N}$.
Then, $D(G)=1+\sum_{i=1}^r(p^{e_i}-1).$ (\cite{Olson})

\item $\eta(C_{n_{1}}\oplus C_{n_{2}})=2n_{1}+n_{2}-2$. (\cite{GZ})

\item Let $G=H\oplus C_n$ with $\exp(H)\mid n\geq 2$. Then, $\eta(G)\geq 2(D(H)-1)+n$. (\cite{affine caps})

\item Every sequence $S$ over $C_{p^n}\oplus H$ of length $|S|=
2p^n+D(H)-2$ contains a zero-sum subsequence $T$ of length $|T|\in
\{p^n,2p^n\}$.

\item $\eta(C_{mp^n}\oplus H)\leq mp^n+p^n+D(H)-2$.

\end{enumerate}
\end{lemma}

\proof \noindent{6.}
\  Let $S=g_1\cdot \ldots \cdot g_{\ell}$ be a sequence over
$G=C_{p^n}\oplus H$ of length $\ell=|S|=2p^n+D(H)-2$. Let $\alpha_i=
{\small\left(\!\!\!\begin{array}{c} 1
\\ g_i
\end{array}\!\!\!\right)}\in C_{p^n}\oplus C_{p^n}\oplus H$ with $1\in C_{p^n}$. Then, $\alpha_1\cdot \ldots \cdot \alpha_{\ell}$ is a sequence over $C_{p^n}\oplus G$ of length
$\ell=p^n+p^n+D(H)-2=D(C_{p^n}\oplus G).$ Therefore, $\alpha_1\cdot
\ldots \cdot \alpha_{\ell}$ contains a nonempty zero-sum subsequence
$W$(say). By the making of $\alpha_i$ we infer that $|W|=p^n$ or
$|W|=2p^n$. Let $T$ be the
 subsequence of $S$ which corresponds to $W$. Then
$T$ is a zero-sum subsequence of $S$ of length $|T|\in \{p^n,
2p^n\}.$

\medskip
\noindent{7.}
\  We first consider the case that $m=1$. Let $G=C_{p^n}\oplus H$.
We want to prove that $\eta(G)\leq 2p^n+D(H)-2$.

Let $S=g_1\cdot \ldots \cdot g_{\ell}$ be a sequence over
$G=C_{p^n}\oplus H$ of length $\ell=|S|=2p^n+D(H)-2$. We need to
show that $S$ contains a short zero-sum subsequence. It follows from
Conclusion 6 that $S$ contains a zero-sum subsequence of length
$|T|\in \{p^n, 2p^n\}$. If $|T|=p^n$ then $T$ itself is a short
zero-sum sequence over $G$ and we are done. Otherwise,
$|T|=2p^n>p^n+D(H)-1=D(G)$. Therefore, $T$ contains a proper
zero-sum subsequence $T'$. Now either $T'$ or $TT'^{-1}$ is a short
zero-sum subsequence of $S$. This proves that $\eta(C_{p^n}\oplus
H)\leq 2p^n+D(H)-2$. By Lemma \ref{Lemma pho(G) and pho(H)}, we have
\begin{align*}
  \eta(C_{p^nm}\oplus H)
  & \leq (\eta(C_{m})-1)\exp(C_{p^n}\oplus H)+\eta(C_{p^n}\oplus H)\\
  & \leq (m-1)p^n+2p^n+D(H)-2\\
  &=mp^n+p^n+D(H)-2.
\end{align*} \qed

\begin{lemma} \label{2exp(G)}  Let $G$ be a  finite abelian group. Then $[D(G)+1, \min \{2\exp(G)+1,
\eta(G)-1\}]\subset C_0(G)$.
\end{lemma}

\proof If $[D(G)+1, \min \{2\exp(G)+1, \eta(G)-1\}] = \emptyset$
then the conclusion of this lemma hold true trivially. Now assume
that $[D(G)+1, \min \{2\exp(G)+1, \eta(G)-1\}] \neq \emptyset$. Let
$S$ be an arbitrary zero-sum sequence over $G$ of length $|S|\in
[D(G)+1,\min\{2\exp(G)+1,\eta(G)-1\}]$. It suffices to show that $S$
contains a short zero-sum subsequence. Since $|S|\geq D(G)+1$, it
follows that $S$ contains a  zero-sum subsequence $T$ of length
$|T|\in [1, |S|-1]$. Then $\sigma(ST^{-1})=0$. Since $|S| \leq
2\exp(G)+1$, we infer that $|T|\in [1,\exp(G)]$ or $|ST^{-1}|\in
[1,\exp(G)]$. This proves the lemma.  \qed

\medskip
\noindent {\sl Proof of Theorem \ref{mainthm2}}.\ \noindent{1.}\ \
By the definition of $C_0(G)$ we have, $C_0(G)\subset [D(G)+1,
\eta(G)-1]$. So, we need to show
$$
[D(G)+1, \eta(G)-1] \subset C_0(G).
$$

Suppose first that
$$G=C_n\oplus C_n.$$ By Lemma \ref{Davenportconstant}, we have
$D(G)=2n-1$ and $\eta(G)=3n-2$. Let $S$ be a zero-sum sequence over
$G$ of length $|S|\in [2n, 3n-3]$. We need to show $S$ contains a
short zero-sum subsequence. We may assume that
$$\mathsf v_0(S)=0.$$ Let $T=S\cdot 0^{3n-2-|S|}$. Then $|T|=3n-2$
and $T$ contains a zero-sum subsequence $T'$ of length $|T'|\in
\{n,2n\}$ by Lemma \ref{nd}. If $|T'|=n$ then $T'0^{-\mathsf
v_0(T')}$ is a short zero-sum subsequence of $S$ and we are done.
So, we may assume that $|T'|=2n$. Let $T''=TT'^{-1}$. Now $T''$ is a
zero-sum subsequence of $T$ of length $|T''|=n-2$. If $T''$ contains
at least one nonzero element then $T''0^{-\mathsf v_0(T'')}$ is a
short zero-sum subsequence of $S$ and we are done. So, we may assume
that $T''=0^{n-2}$. This forces that $T'=S$. It follows from
$D(G)=2n-1$ that $S$ contains a zero-sum subsequence $S_0$ of length
$|S_0|\in [1,2n-1]$. Therefore, either $S_0$ or $SS_0^{-1}$ is a
short zero-sum subsequence of $S$.

Now suppose that $$G=C_n\oplus C_m$$ with $n\mid m$ and $$n<m.$$ By
Lemma \ref{Davenportconstant}, we have that $D(G)=n+m-1<2m$ and
$2m+1>2n+m-2=\eta(G)$. It follows from Lemma \ref{2exp(G)} that
$[D(G)+1, \eta(G)-1]\subset C_0(G)$.

\medskip
\noindent{2.}\ \ By Lemma \ref{Davenportconstant}, we have
that $D(C_{mp^n}\oplus H)=mp^n+D(H)-1$ and $\eta(C_{p^nm}\oplus H)
\leq mp^n+p^n+D(H)-2.$

Suppose $m\geq 2.$ Then $\eta(C_{p^nm}\oplus H)\leq 2mp^n$.
Similarly to the proof of Conclusion 1, we can prove that
$[D(C_{p^nm}\oplus H)+1, \eta(C_{p^nm}\oplus H)-1]\subset C_{0}(G)$,
and we are done. So, we may assume $$m=1.$$ Then $\eta(C_{p^n}\oplus
H)\leq 2p^n+D(H)-2$ and the proof is similar to that of 1 by using
Conclusion 6 of Lemma \ref{Davenportconstant}.

\medskip
\noindent{3.}\ \
It is just Proposition \ref{propaa}.

\medskip
\noindent{4.}\ \ It is easy to verify that
$C_0(C_2^r)\supset\{2^r-3,2^r-2\}=\{\eta(C_2^r)-3,\eta(C_2^r)-2\}$.
So, $C_0(C_2^r)=\{2^r-3,2^r-2\}=\{\eta(C_2^r)-3,\eta(C_2^r)-2\}$
follows from Proposition \ref{zz}. \qed

\section{Proof of Theorem \ref{mainthm1}}

We  need the following result which states that Property C is
multiple.

\begin{lemma} (\cite{GGS07} Theorem 3.2) \label{propertyC-multiple}  Let
$G=C_{mn}^r$ with $m,n,r \in \mathbb{N}$. If both $C_m^r$ and
$C_n^r$ have Property C and
$$
\frac{\eta(C_m^r)-1}{m-1}=\frac{\eta(C_n^r)-1}{n-1}=\frac{\eta(C_{mn}^r)-1}{mn-1}=c
$$
for some $c\in \mathbb{N}$ then $G$ has Property C.
\end{lemma}

\begin{lemma} \label{s-multiple}  If $\frac{\eta(C_m^r)-1}{m-1}=\frac{\eta(C_n^r)-1}{n-1}=c$ for some $c\in \mathbb{N}$ and if $\eta(C_{mn}^r)\geq
c(mn-1)+1$ then $\eta(C_{mn}^r)=c(mn-1)+1.$
\end{lemma}

\proof The lemma follows from Lemma \ref {Lemma pho(G) and pho(H)}.
\qed

\medskip
We also need the following easy lemma
\begin{lemma} \label{s(G)-eat(G)} (\cite{GR} Lemma 4.2.2) Let $G$ be a finite abelian group.
Then, $\mathsf s(G)\geq \eta(G)+\exp(G)-1.$
\end{lemma}

\begin{proposition} \label{Theorem General}  Let $n=3m$, where $m$ is an odd positive
integer. Then,
\begin{enumerate}
\item If $\eta(C_m^3)=8m-7$ then $\eta(C_n^3)-2\in C_0(C_n^3)$.

\item If $\eta(C_m^4)=19m-18$ then
$\{\eta(C_n^4)-2,\eta(C_n^4)-1\}\subset C_0(C_n^4)$.
\end{enumerate}
\end{proposition}

\proof \noindent{1.}\ \  From Lemma \ref{Lemma
pho(G) and pho(H)} and Lemma \ref{Lemma list of contants} we deduce
that
$\frac{\eta(C_n^3)-1}{n-1}=\frac{\eta(C_m^3)-1}{m-1}=\frac{\eta(C_3^3)-1}{3-1}=8$
. Now  $\eta(C_n^3)-2\in C_0(C_n^3)$ follows from Lemma \ref{Lemma
list of contants}, Lemma \ref{Lemma C33-2 and C34-2} and Lemma
\ref{Theorem general for C_0} with $G_2=C_3^3$ and $t_1=t_2=2$.

\medskip
\noindent{2.} \ \  The proof is similar to that of Conclusion 1.
\qed

\begin{proposition} \label{Theorem C35^3-2 in}  Let $\alpha, \beta \in \mathbb{N}_0$ with $\alpha\geq 1$.
Then,
\begin{enumerate}
\item If $\alpha+\beta\geq 2$ then $\{\eta(C_{3^{\alpha}
5^{\beta}}^3)-2,\eta(C_{3^{\alpha} 5^{\beta}}^3)-1\}\subset
C_0(C_{3^{\alpha} 5^{\beta}}^3)$.

\item $\{\eta(C_{3^{\alpha}}^4)-2,\eta(C_{3^{\alpha}}^4)-1\}\subset
C_0(C_{3^{\alpha}}^4)$.
\end{enumerate}
\end{proposition}

\proof The proposition follows from Lemma \ref{Lemma list of
contants}, Lemma \ref{s-multiple}, Lemma \ref{Prog} and Proposition
\ref{Theorem General}. \qed

\medskip
\noindent We shall show the following property which is easier to be
checked than Property C can also be used to study $C_0(G)$.

\medskip
\noindent{\bf Property $D_0$:}   {\sl Let $c, n\in \mathbb{N}$. We
say that $C_n^r$ has property $D_0$ with respect to $c$ if every
sequence of the form $g\prod_{i=1}^{c}g_i^{n-1}$ contains a zero-sum
subsequence of length exactly $n$, where $g, g_i\in C_n^r$ for all
$i\in [1, c]$}.

\begin{lemma}\label{EGZ1}(\cite{FGZ}, page 8)  Let $m=3^a5^b$ with $a,b$ nonnegative integers. Let $n\geq 65$
be an odd positive integer such that $C_p^3$ has Property $D_0$ with
respect to 9 for all prime divisors $p$ of $n$. If $$m \geq
\frac{2\times 5^7n^{17}}{(n^2-7)n-64}$$ then $\mathsf
s(C_{mn}^3)=9mn-8.$
\end{lemma}

\medskip
\begin{proposition}\label{22}   Let $m=3^{\alpha}5^{\beta}$ with $\alpha\in \mathbb{N}$ and $\beta\in \mathbb{N}_{0}$. Let $n\geq 65$ be an
odd positive integer such that $C_{p}^{3}$ has Property $D_{0}$ with
respect to 9 for all prime divisors $p$ of $n$. If $$m\geq
\frac{6\times 5^7n^{17}}{(n^2-7)n-64}+3$$ then
$\eta(C_{mn}^{3})-2\in C_{0}(C_{mn}^{3}).$
\end{proposition}

\proof Let $m'=\frac{m}{3}$. Then
$m'=3^{\alpha-1}5^{\beta}\geq\frac{2\times 5^7n^{17}}{(n^2-7)n-64}$
and $\alpha-1\geq 0$.

By Lemma \ref{EGZ1} and Lemma \ref{s(G)-eat(G)} we have
$s(C_{m'n}^{3})=9m'n-8$ and $\eta(C_{m'n}^{3})\leq 8m'n-7$. It
follows from Lemma \ref{Lemma list of contants} that
$\eta(C_{m'n}^{3})=8m'n-7$. Since $\eta(C_{3}^{3})=8\times 3-7$ and
$\eta(C_{m'n}^{3})=8m'n-7$, it follows from Lemma \ref{s-multiple}
that $\eta(C_{mn}^{3})=8mn-7$. What's more, $C_{3}^{3}$ has Property
$C$ and $\eta(C_{3}^{3})-2\in C_{0}(C_{3}^{3})$ by Lemma \ref{Lemma
C33-2 and C34-2}. Therefore,  $\eta(C_{mn}^{3})-2\in
C_{0}(C_{mn}^{3})$ by Lemma \ref{Theorem general for C_0}.\qed

\medskip
{\sl Proof of   Theorem \ref{mainthm1}. }\noindent{1.} It follows
from Theorem \ref{mainthm2} that $C_0(G)=[D(G)+1,
\eta(G)-1]=[n_2+n_1, n_2+2n_1-3]\neq \emptyset$.

\medskip
\noindent{2.} If
$a\geq 1$ then it follows from Proposition \ref{Theorem C35^3-2 in}
and Lemma \ref{Lemma C33-2 and C34-2}. Now assume $b\geq 2$. Since
$\eta(C_{3^a5^b}^{3})=8(3^a5^b-1)+1$, it follows from Lemma
\ref{Prog} that $\eta(C_{3^a5^b}^{3})-1\in C_{0}(C_{3^a5^b}^{3})$.

\medskip
\noindent{3.} Let $G_{1}=C_{3\times 2^{a-3}}^{3}$ and
$G_{2}=C_{8}^{3}$. Then $\eta(G_{1})=7(3\times 2^{a-3}-1)+1$,
$\eta(G_{2})=7\times (8-1)+1$ and $G_{2}$ has Property C. Therefore,
$\eta(C_{8}^{3})-1\in C_{0}(C_{8}^{3})$ by Lemma \ref{Prog}. So,
$\eta(C_{3\times 2^a}^{3})-1\in C_{0}(C_{3\times 2^a}^{3})$ by Lemma
\ref{Theorem general for C_0}.

\medskip
\noindent{4.} The result follows from Proposition \ref{Theorem
C35^3-2 in}.

\medskip
\noindent{5.} Let $G=C_{2^a}^{r}$ with $3\leq r\leq a$.
 Since $\eta(C_{2^a}^{r})=(2^r-1)(2^a-1)+1$, $2^r-1<2^a$ and $C_{2^a}^{r}$ has Property C, we have that $\eta(C_{2^a}^{r})-1\in C_{0}(C_{2^a}^{r})$ by Lemma \ref{Prog}.

If $G=C_{2}^{r}$ and $r\geq 3$, then it follows from Theorem
\ref{mainthm2}.

\medskip
\noindent{6.} Let $m=3^{n_1}5^{n_2}$ and
$n=7^{n_3}11^{n_4}13^{n_5}$. It follows from $n_3+n_4+n_5\geq 3$
that $n>65$. By the hypothesis of $n_1+n_2\geq 11+34(n_3+n_4+n_5)$
we infer that, $m=3^{n_1}5^{n_2}\geq 3^{n_1+n_2}\geq
3^{11}3^{34(n_3+n_4+n_5)}>4\times 5^8\times
13^{14(n_3+n_4+n_5)}>4\times 5^8n^{14}>\frac{6\times 5^7
n^{17}}{(n^2-7)n-64}+3$. Since it has been proved that every prime
$p\in \{3, 5, 7, 11, 13\}$ has Property $D_{0}$ with respect to 9 in
\cite{FGZ}, it follows from Proposition \ref{22}  that
$\eta(C_{k}^{3})-2\in C_{0}(C_{k}^{3})$.

\medskip
\noindent{7.} It follows from Conclusion 2 of Theorem \ref{mainthm2}
that $C_0(G)=[D(G)+1, \eta(G)-1]$. By 3 and 5 of Lemma
\ref{Davenportconstant} we have that, $\eta(G)-1\geq
mp^n+2(D(H)-1)-1\geq mp^n+D(H)=D(G)+1$. Therefore, $C_0(G)\neq
\emptyset$. \qed

\section{ Concluding Remarks and Open Problems}

\begin{proposition}\label{C0(G) contains no consecutive}  Let $G$ be a non-cyclic finite abelian group with $\exp(G)=n$. Then
$C_0(G)\cup \{\eta(G)\}$ doesn't contain $n+1$ consecutive integers.
\end{proposition}

\proof Assume to contrary that $[t,t+n]\subset C_0(G)\cup
\{\eta(G)\}$ for some $t\in \mathbb{N}$. By the definition of
$C_0(G)$ we have that $t+n-1<\eta(G)$. So, we can choose a short
free sequence $T$ over $G$ of length $|T|=t+n-1$. It follows from
$t+n-1\in C_0(G)\cup \{\eta(G)\}$ that $\sigma(T)\neq 0$. Let
$g=\sigma(T)$ and let $S=T\cdot (-g)$. Since $|S|=t+n\in C_0(G)\cup
\{\eta(G)\}$, $S$ contains a short zero-sum subsequence $U$ with
$(-g)\mid U$. Note that $t\leq |S\cdot U^{-1}|\leq t+n-2$ and
$\sigma(S\cdot U^{-1})=0$. It follows from $[t,t+n]\subset
C_0(G)\cup \{\eta(G)\}$ that $S\cdot U^{-1}$ contains a short
zero-sum subsequence, which is a contradiction with $S\cdot
U^{-1}\mid T$. \qed

\medskip
Proposition \ref{C0(G) contains no consecutive} just asserts that
$C_0(G)$ can't contain any interval of length more than $\exp(G)$.
Proposition \ref{zz} shows that $C_0(C_n^r)$ could not contain
integers much smaller than $\eta(C_n^r)-1$. So, it seems plausible
to suggest

\begin{conjecture} \label{Conjecture general}  Let $G\neq C_{2}\oplus C_{2m}, m\in \mathbb{N}$ be a finite non-cyclic abelian group. Then

$C_0(G)\subset [\eta(G)-(\exp(G)+1),\eta(G)-1]$.
\end{conjecture}

Conjecture \ref{Conjecture general} and Conjecture \ref{conjecture1}
suggest the following
\begin{conjecture} \label{Conjecture genera2}  Let $G\neq C_{2}\oplus C_{2m}, m\in \mathbb{N}$ be a finite non-cyclic abelian group.
Then $1\leq |C_0(G)|\leq \exp(G).$
\end{conjecture}

\begin{conjecture} \label{Conjecture genera3}
 $C_0(G)=[\min C_0(G), \max C_0(G)]$.
\end{conjecture}

The following notation concerning the inverse problem on $\mathsf
s(G)$ was introduced in \cite{GG06}.

\noindent{\bf Property D:}  {\sl We say the group $C_{n}^{r}$ has
property D if $\mathsf s(C_{n}^r)=c(n-1)+1$ for some positive
integer
 $c$, and every  sequence $S$ over $C_n^r$ of length $|S|=c(n-1)$ which contains no zero-sum subsequence of length $n$  has the
 form $S=\prod_{i=1}^{c}g_{i}^{n-1}$ where $g_1,\ldots,g_c$ are pairwise distinct elements of $C_n^r$.}

 \medskip
\begin{conjecture} (\cite{GG06}, Conjecture 7.2) \label{Conjecture genera4}
 Every group $C_n^r$ has Property D.
\end{conjecture}

It has been proved (\cite{GG06}, Section 7) that Conjecture
\ref{Conjecture genera4}, if true would imply
\begin{conjecture} \label{Conjecture genera5}
 Every group $C_n^r$ has Property C.
\end{conjecture}

Suppose that Conjecture \ref{Conjecture genera5} holds true for all
groups of the form $C_n^r$. For fixed $n,r\in \mathbb{N}$ and any
$a\in \mathbb{N}$ we have that $\eta(C_{n^a}^r)=c(n^a,r)(n-1)+1$,
where $c(n^a,r)\in \mathbb{N}$ depends on $n^a$ and $r$. By Lemma
\ref{Lemma pho(G) and pho(H)} we obtain that the sequence
$\{c(n^a,r)\}_{a=1}^{\infty}$ is decreasing. Therefore,
$c(n^a,r)\leq n^a$ for all sufficiently large $a$. Hence, by Lemma
\ref{Prog} we infer that $\eta(C_{n^a}^r)-1\in C_0(C_{n^a}^r)$ for
all sufficiently large $a\in \mathbb{N}$.

\bigskip
\noindent {\bf Acknowledgments}. This work has been supported by the
PCSIRT Project of the Ministry of Science and Technology, the
National Science Foundation of China with grant no. 10971108 and
11001035, and the Fundamental Research Funds for the Central
Universities.


\bigskip

\begin{thebibliography}{99}

\bibitem{Davenport}
G. Bhowmik and J.C. Schlage-Puchta, Davenport' s constant for groups
of the form $\mathbb{Z}_3\oplus \mathbb{Z}_3 \oplus\mathbb{Z}_{3d}$,
\emph{Additive Combinatorics} \textbf{43} (2007), 307--326.

\bibitem{Lower bounds}
C. Elsholtz, Lower bounds for multidimentional zero sums,
\emph{Combinatorica} \textbf{24} (2004), 351--358.

\bibitem{Restricted size IV}
R. Chi, S.Y. Ding, W.D. Gao, A.~Geroldinger and W.A. Schmid, On
zero-sum subsequences of restricted size IV, \emph{Acta Math.
Hungar.} \textbf{107} (2005), 337--344.

\bibitem{Game}
B.L. Davis, D. Maclagan and R. Vakil, The card game set, \emph{Math.
Intelligencer} \textbf{25} (2003), 33--40.

\bibitem{Weight rules}
Y. Edel, Sequences in abelian groups G of odd order without zero-sum
subsequences of length exp(G), \emph{Des. Codes Cryptogr.}
\textbf{47} (2008), 125--134.

\bibitem{affine caps}
Y. Edel, C. Elsholtz, A. Geroldinger, S. Kubertin and L. Rackham,
Zero-sum problems in finite abelian groups and affine caps, \emph{Q.
J. Math.} \textbf{58} (2007), 159--186.

\bibitem{EK69}
P. van Emde Boas and  D. Kruyswijk, A combinatorial problem on
finite abelian groups, II. \emph{Math. Centre Report} ZW-1969-008.

\bibitem{FGZ}
Y.S. Fan, W.D. Gao and Q.H. Zhong, On the Erdos-Ginzburg-Ziv
constant of finite abelian groups of high rank, \emph{J. Number
Theory} \textbf{131} (2011), 1864--1874.

\bibitem{Gao97}
W.D. Gao, On a combinatorial problem connected with factorizations,
\emph{Colloq. Math.} \textbf{72} (1997), 251--268.

\bibitem{GG06}
W.D. Gao and A. Geroldinger, Zero-sum problems in finite abelian
groups{\rm \,:} a survey, \emph{Expo. Math.} \textbf{24} (2006),
337--369.

\bibitem{GHST}
W.D. Gao, Q.H. Hou, W.A. Schmid and R. Thangadurai, On short
zero-sum subsequences II, \emph{Integers} \textbf{7} (2007), A21.

\bibitem{GGS07}
W.D. Gao, A. Geroldinger and W.A. Schmid, Inverse zero-sum problems,
\emph{Acta Arith.} \textbf{128} (2007), 245--279.

\bibitem{GH06}
A. Geroldinger and F. Halter-Koch, Non-{U}nique {F}actorizations.
{A}lgebraic, {C}ombinatorial and {A}nalytic {T}heory, \emph{Pure and
Applied Mathematics}, vol. 278, Chapman \& Hall/CRC, (2006).

\bibitem{GZ}
A. Geroldinger, Additive group theory and non-unique factorizations,
\emph{Combinatorial Number Theory and Additive Group Theory} (2009),
1--86.

\bibitem{GGS}
A. Geroldinger, D.J. Grynkiewicz and  W.A. Schmid, Zero-sum problems
with congruence conditions, \emph{Acta Math. Hungar.} \textbf{131}
(2011), 323--345.

\bibitem{GR}
A. Geroldinger and I.Z. Ruzsa, Combinatorial number theory and
additive group theory, Birkhauser (2009).


\bibitem{Harborth I}
H. Harborth, Ein Extremalproblem f\"{u}r Gitterpunkte, \emph{J.
Reine Angew. Math.} \textbf{262} (1973), 356--360.

\bibitem{Harborth}
H. Harborth, Ein Extremalproblem f\"{u}r Gitterpunkte, Ph.D. Thesis,
Technische Univesit\"{a}t Braunschweig (1982).

\bibitem{Kemnitz}
A. Kemnitz, On a lattice point problem, \emph{Ars Combin.} 16b
(1983), 151--160.

\bibitem{Olson}
J.E. Olson, A combinatorial problem on finite abelian groups, I; II,
\emph{J. Number Theory} \textbf{1} (1969), 8--10; 195--199.

\end{thebibliography}
\end{document}